\documentclass{amsart}

\usepackage{amssymb,latexsym}

\newtheorem{theorem}{Theorem}
\newcommand{\bt}{\begin{theorem}}
\newcommand{\et}{\end{theorem}}
\newtheorem{lemma}{Lemma}
\newcommand{\bl}{\begin{lemma}}
\newcommand{\el}{\end{lemma}}
\newtheorem{corollary}{Corollary}
\newcommand{\bc}{\begin{corollary}}
\newcommand{\ec}{\end{corollary}}
\newtheorem{problem}{Problem}
\newcommand{\bprob}{\begin{problem}}
\newcommand{\eprob}{\end{problem}}
\newcommand{\beq}{\begin{equation}}
\newcommand{\eeq}{\end{equation}}
\newcommand{\benum}{\begin{enumerate}}
\newcommand{\eenum}{\end{enumerate}}
\newcommand{\N}{\ensuremath{ \mathbf N }}
\newcommand{\Z}{\ensuremath{\mathbf Z}}

\newcommand{\PP}{\ensuremath{ \mathbf P}}

\newcommand{\mci}{\ensuremath{ \mathcal I}}

\newcommand{\bmat}{\left(\begin{matrix}}
\newcommand{\emat}{\end{matrix}\right)}
\newcommand{\bsmallmat}{\left(\begin{smallmatrix}}
\newcommand{\esmallmat}{\end{smallmatrix}\right)}

\DeclareMathOperator{\qqand}{\qquad\text{and}\qquad}

\title[Dyson transform]{Additive number theory and the Dyson transform}
\author{Melvyn B.   Nathanson}
\address{Lehman College (CUNY), Bronx, NY 10468}
\email{melvyn.nathanson@lehman.cuny.edu}

\date{\today}



 \thanks{MBN supported in part by  PSC-CUNY Research Award Program grant 66197-00 54.}
\begin{document}

\begin{abstract}
In 1942 Mann solved a famous problem, the $\alpha+\beta$ conjecture, 
about the lower bound of the Shnirel'man density of sums of sets of positive integers.  
In 1945, Dyson generalized Mann's theorem and obtained a lower bound for the 
Shnirel'man density of rank $r$ sumsets.  His proof introduced the Dyson transform, 
an important tool in additive number theory.
This paper explains the background of Dyson's work, gives Dyson's proof of his theorem, 
and includes several applications of the Dyson transform, 
such as  Kneser's inequality for sums of finite subsets of an arbitrary 
additive abelian group. 
\end{abstract}

\maketitle

\section{A theorem on the densities of sets of integers}

In 1945 Freeman Dyson~\cite{dyso45} published a paper in additive number 
theory that gave a lower bound for the Shnirel'man density of $r$-rank 
sums of sets of  nonnegative integers. 
In a commentary on this paper in his \emph{Selected Works}~\cite[pp. 5--6]{dyso96}, 
Dyson wrote an autobiographical history of this result.  
\begin{quotation}
This paper was inspired by \ldots the little book, 
``\" Uber einige neuere Fortschritte der additiven Zahlentheorie,'' 
by Edmund Landau, [Cambridge University Press, 1937]\nocite{land37}.  
As soon as I arrived in Cambridge in 1941, I scoured the book-shops for bargains 
and bought Landau's book 
for six shillings and sixpence\ldots. 
His Chapter 4 describes the ``alpha-beta conjecture'', which was in 1937 
a famous unsolved problem.  
The conjecture concerns the densities of infinite sets of positive integers\ldots.  
I laid the problem in my heart and worked hard to find a proof\ldots.

While I was at Cambridge, I learned that Henry Mann \ldots had proved 
the alpha-beta conjecture, [Annals of Mathematics, \textbf{43}, 523--527 (1942)]\ldots.  
I was sorry that Mann had beaten me in the race to find a proof, 
but I was glad to see that his strengthening of the conjecture still left me something to do\ldots.  
The generalization of Mann's theorem to the sum of more than two sets 
remained to be proved\ldots. 
I took the problem with me when I left Cambridge in 1943 to work in the 
Operational Research Section of the Royal Air Force Bomber Command 
at High Wycombe.  I solved it in the early months of 1944 during a particularly grim 
period in the history of the Bomber Command\ldots.  
The Mann Theorem helped to keep me sane amid the insanities of the bombing campaign.
\end{quotation}

Dyson's proof of the generalization of Mann's theorem was inductive, and 
involved the replacement of one pair $(A_{\ell},A_h)$ of  sets of positive integers with 
another pair $(A_{\ell}',A_h')$ of sets such that 
$|A_h' \cap \{1,\ldots, g\} | < |A_h\cap  \{1,\ldots, g\} |$ 
for some positive integer $g$.  
Through the work of Martin Kneser, a variant of this transformation, 
variously called the \emph{Dyson transform}, 
the \emph{$e$-transform}, and the \emph{Dyson $e$-transform},  
has become a fundamental tool in additive number theory.  
This paper explains the mathematical context of Dyson's paper, 
presents a proof of Dyson's addition theorem, and describes some 
applications of the Dyson transform.  
For a history of the $\alpha + \beta$ conjecture, see Erd\H os and Niven~\cite{erdo-nive46}.

\section{Shnirel'man density and sums of primes}

Associated to every set $A$ of integers  is the number $\sigma(A)$, called the 
\emph{Shnirel'man density} of $A$.  
The \emph{Dyson transform}\index{Dyson transform} was introduced to obtain 
a lower bound for Shnirel'man density.  To appreciate Dyson's work, we 
need to understand  the importance of Shnirel'man density in additive number theory. 

Let $(A_1, A_2, \ldots, A_h)$ be an $h$-tuple of sets of positive integers.
The \emph{sumset} 
\[
\sum_{i=1}^h A_i = A_1 + A_2 + \cdots + A_h 
\]
is the set of all integers 
that can be written in the form $a_1+a_2+\cdots + a_h$, where $a_i \in A_i \cup \{0\}$ 
for all $i \in \{1,2,\ldots, h\}$.  The \emph{$h$-fold sumset} of a set $A$ is the set 
\[
hA = \underbrace{A + A + \cdots + A}_{\text{$h$ summands}}.
\]
Let $\N_0 = \{0,1,2,\ldots\}$ be the set of nonnegative integers.  
The set $A$ of positive integers is a \emph{basis of order $h$}\index{basis} if 
every nonnegative integer is a sum of $h$ elements of $A \cup \{0\}$, 
or, equivalently, if $hA = \N_0$.  
The set $A$   is a \emph{basis} if $A$ is a basis of order $h$ for 
some positive integer $h$. 
For example, by Lagrange's theorem, the set  of squares is a basis of order 4. 
Waring's problem states that, for every integer $k \geq 2$, there is a positive 
integer $g(k)$ such that the set of positive $k$th powers is a basis of order $g(k)$.

Let $\PP = \{2,3,5,7,\ldots\}$ be the set of prime numbers and let 
$\widehat{\PP} = \{1\} \cup \PP$.  
The Goldbach conjecture that every even integer greater than 2 is the sum of two primes  
implies that the set $\widehat{\PP} $ is a basis of order 3.  
In 1930, in an extraordinarily original and powerful work, the Soviet mathematician 
Lev Genrikhovich Shnirel'man~\cite{shni30,shni33,land30} 
proved that  there exists an integer $h$ such that every integer greater than 1 
is the sum of at most $h$ prime numbers. 
Equivalently, the set $\widehat{\PP}$ is a basis.  
To obtain this result, Shnirel'man 
introduced a new density, which we call \emph{Shnirel'man density}\index{Shnirel'man density}, 
on sets of integers and proved that every set of positive integers 
with positive Shnirel'man density is a basis.  
He used the Brun sieve to prove that the sumset $h'\widehat{\PP}$ 
has positive Shnirel'man density for some positive integer $h'$, and, 
\emph{mirabile dictu}, obtained his result.  

Recent results on the Goldbach conjecture include Ramar\' e's proof~\cite{rama95} that 
every even integer $n \geq 2$ is the sum of at most 6 primes, Tao's proof~\cite{tao14b}
that every odd integer $n \geq 3$ is the sum of at most 5 primes, 
and Helfgott's proof~\cite{helf15b} that every odd integer $n \geq 7$ is the sum of 3 primes.

To define Shnirel'man density, we use the 
 \emph{counting function}  $A(x)$ of a set $A$ of integers.
For all real numbers $x$, this function counts the number of positive integers 
in $A$ that do not exceed $x$, that is,
\[
A(x) = \sum_{\substack{a \in A \\ 1 \leq a \leq x}} 1.
\]
Note that $A(x) = 0$ if $x < 1$. 
The \emph{Shnirel'man density}\index{Shnirel'man density} 
of the set $A$ is
\[
\sigma(A) = \inf_{n\in \N} \frac{A(n)}{n}.
\]

Shnirel'man density has some usual and some unusual properties.  
Let $A$ be set of positive integers. 
We have $0 \leq \sigma(A) \leq 1$.   
The set of odd integers has Shnirel'man density $1/2$.  
However,  $\sigma(A) = 0$ if $1 \notin A$ (because $A(1) = 0$) 
and so  the set of even integers has Shnirel'man density $0$.  
A fundamental fact about Shnirel'man 
density is that $\sigma(A) = 1$ if and only if $A$ contains all positive integers. 
Thus, the set $A$ is basis of order $h$ for some positive integer $h$ 
if and only if $\sigma(hA) = 1$.  
An important result is Shnirel'man's theorem that a set $A$ with $\sigma(A) > 0$  
is a basis of order $h$ for some $h$ .  
The proof begins with a simple counting argument. 
If $\sigma(A) = \alpha \geq 1/2$, then, for every positive integer $n$, 
we have 
\[
A(n) \geq \alpha n \geq \frac{n}{2}.
\]
Let 
\[
1 \leq a_1 < \cdots < a_k \leq n
\]
be the $k = A(n)$ positive integers in $A$ that do not exceed $n$.
If $n \in A$, then $n = 0+n \in 2A$.  If $n \notin A$, then $a_k \leq n-1$ and so 
\[
1 \leq a_1 < \cdots < a_k \leq n-1  
\]
and  
\[
1 \leq n-a_k < \cdots < n-a_1 \leq n-1.
\]
Thus, the set $\{1,\ldots, n-1\}$ contains the subsets $A_1 = \{a_i:i=1,\ldots, k\}$ and 
$A_2 = \{n-a_j:j=1,\ldots, k\}$.  
Both sets have cardinality $k = A(n) \geq n/2$, and so they cannot be disjoint.  
Thus, there exist $i,j \in \{1,\ldots, k\}$ such that 
\[
a_i = n-a_j \in A_1 \cap A_2
\]
and $n=a_i+a_j \in 2A$.
Thus, $\sigma(A) \geq 1/2$ implies $\sigma(2A) = 1$ and so $A$ is a basis of order 2.  

The crux of Shnirel'man's proof is a lower bound for the  density of an $h$-fold sumset.
Let  $A$ and $B$ be sets of positive integers.  
If $\sigma(A) = \alpha$ and $\sigma(B) = \beta$, then 
\beq            \label{Dyson:ShnirelmanIneq}
\sigma(A+B) \geq \alpha + \beta- \alpha \beta.  
\eeq
(For a proof of this inequality, see Nathanson~\cite[Section 7.4]{nath96aa}.)
Equivalently, 
\[
1 - \sigma(A+B) \leq (1-\alpha)(1-\beta).
\]
It follows by induction that if $A_1,\ldots, A_h$ are sets of integers 
with $\sigma(A_i) = \alpha_i$ for all $i = 1,\ldots, h$, then 
\[
1 - \sigma\left(A_1 + \cdots + A_h \right) \leq \prod_{i=1}^h (1-\alpha_i).
\]
Let $A$ be a set of integers with $\sigma(A) = \alpha$.  
Applying~\eqref{Dyson:ShnirelmanIneq} with $A_i = A$ 
and $\alpha_i = \alpha$ for all $i \in \{1,\ldots, h\}$, 
we obtain 
\[
1 - \sigma\left(hA \right) \leq  (1-\alpha)^h.  
\]
Equivalently,
\[
\sigma\left(hA \right) \geq  1 - (1-\alpha)^h.
\]
If $\alpha > 0$, then $\lim_{h\rightarrow \infty} (1-\alpha)^h = 0$ 
and so $\sigma(hA) \geq 1/2$ for all   $h \geq h_0$.   
It follows that $\sigma(2hA)  = 1$ for all $h \geq 2h_0$.  Thus, to prove that a set $A$ 
is a basis, it suffices to prove that the sumset $hA$ has positive
Shnirel'man density for some positive integer $h$. 

Let $A$ and $B$ be sets of positive integers with $\sigma(A) = \alpha$ 
and $\sigma(B) = \beta$.    
Khinchin~\cite{khin32,khin33} conjectured that Shnirel'man's  density 
inequality~\eqref{Dyson:ShnirelmanIneq} can be strengthened as follows:  
\beq                                \label{Dyson:ShnirelmanIneq-strong}
\sigma(A+B)  \geq \min(1,\alpha+ \beta).
\eeq 

Following work of Besicovich~\cite{besi35}, Brauer~\cite{brau41}, Khinchin~\cite{khin32}, 
Landau~\cite{land30}, and Schur~\cite{schu36}, 
Henry B. Mann~\cite{mann42}  proved a finite inequality that implied Khinchin's conjecture in 1942.  
Let $A+B=C$.  Mann wrote,  
``Stripped of its transcendental content~\eqref{Dyson:ShnirelmanIneq-strong} states that
\[
\frac{C(n)}{n} = 1 \quad \text{or} \quad 
  \geq \min  \frac{A(\ell)}{\ell} + \frac{B(m)}{m}, \qquad 1 \leq m \leq n, \quad 1 \leq \ell \leq n.
\]
We propose to prove the following sharper theorem.''

\bt[Mann~\cite{mann42}]           \label{Dyson:theorem:Mann-FundamentalTheorem}
Let $(A,B)$ be a pair of sets of positive integers and let $A+B = C$.  Then 
\[
\frac{C(n)}{n} = 1
\qquad\text{or} \qquad 
\frac{C(n)}{n}  \geq \min_{\substack{1 \leq m \leq n\\m \notin C}} \frac{A(m)+B(m)}{m}.
\]
\et

The proof of this sharper result, which Mann called the ``Fundamental Theorem,'' 
is elementary but difficult.  
Expositions of this work are in Gelfond-Linnik~\cite{gelf-linn65},   
Halberstam-Roth~\cite{halb-roth66}, and Khinchin~\cite{khin52}.

\bc          \label{Dyson:theorem:Mann-finite}
Let $(A,B)$ be a pair of sets of positive integers and let $A+B = C$.  
If 
\[
A(m) + B(m) \geq \gamma m 
\]
for all $m = 1,\ldots, n$, then 
\[
C(n) \geq \min(1,\gamma) n.
\]
\ec

\begin{proof}
Theorem~\ref{Dyson:theorem:Mann-FundamentalTheorem} 
implies that either $C(n)/n = 1$ or $C(n)/n \geq \gamma$, and so 
$C(n) \geq \min(1,\gamma) n$.  
\end{proof} 

We immediately obtain the $\alpha+\beta$ conjecture. 

\bt[Mann~\cite{mann42}]                     \label{Dyson:theorem:alpha+beta}
Let $(A,B)$ be a pair of sets of positive integers and let $C = A+B$.  
If $\sigma(A) = \alpha$  and  $\sigma(B) = \beta$, then  
\[
\sigma(A+B) \geq \min(1,\alpha+\beta). 
\]
\et

\begin{proof}
If $\sigma(A) = \alpha$ and $\sigma(B) = \beta$, then,  
for all $m = 1,2,\ldots$,  we have 
$A(m) \geq \alpha m$ and $B(M) \geq \beta m$  
and so $A(m) + B(m) \geq (\alpha+\beta) m$.
\end{proof} 

In 1943, Emil Artin and Peter Scherk~\cite{arti-sche43} 
improved Mann's theorem.  
They wrote, ``We found that his method can be simplified considerably 
and even yields some stronger results.'' 
In 1945, Freeman Dyson~\cite{dyso45} refined the Artin-Scherk method and 
extended Mann's $\alpha+\beta$ theorem to $r$-fold 
sums of $h$ sets of positive integers.  

Let $(A_1,\ldots, A_h)$ be an $h$-tuple of sets of positive integers.  
The sets $A_i$ are not necessarily distinct. 
Let $\mci(r,h)$ be the set of all subsets $I$ of $\{1,\ldots, h\}$ of size $r$.  
Note that $\mci(r,h-1) \subseteq \mci(r,h)$ for all $h \geq 1$.  
A \emph{sum of rank $r$}\index{sum of rank $r$} of the $h$-tuple 
$(A_1,\ldots, A_h)$ is a sumset of the form 
\[
S_I = A_{i_1}+\cdots + A_{i_r}
\]
where  $I = \{i_1,\ldots, i_r \} \in \mci(r,h)$.
There are $\binom{h}{r}$ subsets $I$ of $\{1,\ldots, h\}$ of size $r$, 
and so there are 
$\binom{h}{r}$ sums of rank $r$ of the $h$ sets $A_1,\ldots, A_h$.   
These sumsets are not necessarily distinct. 
For every positive integer  $m$,  the \emph{rank $r$ counting function} 
of the  $h$-tuple  $(A_1,\ldots, A_h)$ is 
\[
\phi_r(n) = \sum_{I\in \mci(r,h)} S_I(n) 
\] 
where the sum is over all sums of rank $r$ of the sets $A_1,\ldots, A_h$. 
Note that  
\[
\phi_1(n) = A_1(n)+A_2(n)+\cdots + A_h(n)
\]
and
\[
\phi_h(n) = (A_1+A_2+\cdots + A_h)(n). 
\]

Nathanson proved that Mann's $\alpha+\beta$ theorem 
(Theorem~\ref{Dyson:theorem:alpha+beta}) extends to the 
Shnirel'man density of $r$-fold sumsets of infinite sets.

\bt[Nathanson~\cite{nath1990-78}]       \label{Dyson:theorem:Nathanson}
Let $(A_1,\ldots, A_h)$ be an $h$-tuple of sets of positive integers and let 
$\sigma(A_i) = \alpha_i$ for all $i \in \{1,\ldots,h\}$.  Then
\[
\sum_{I\in \mci(r,h)} \sigma(S_I) \geq  \binom{h-1}{r-1} \min(1,\alpha_1+\cdots + \alpha_h).
\]
\et

Cheo~\cite{cheo52}, Heged{\" u}s-Piroska-Ruzsa~\cite{hege-piro-ruzs98}, 
Lepson~\cite{leps50}, and Nathanson~\cite{nath1990-78} constructed pairs of sets 
$(A,B)$ that show that the lower bound
in Theorem~\ref{Dyson:theorem:alpha+beta} is best possible.  
Nathanson~\cite{nath1990-78} constructed $h$-tuples of sets $(A_1,\ldots, A_h)$ 
that prove that the lower bounds
in Theorem~\ref{Dyson:theorem:Nathanson} are best possible.  

Dyson extended Mann's more difficult Theorem~\ref{Dyson:theorem:Mann-finite} 
to $r$-fold sumsets of finite sets. 

\bt[Dyson~\cite{dyso45}]            \label{Dyson:theorem:Dyson} 
Let $(A_1,\ldots, A_h)$ be an $h$-tuple of sets of  positive integers 
and let $n$ be a positive integer.  If 
\[
\phi_1(m) \geq \gamma m
\]
for all $m \in \{1,\ldots, n\}$ and 
\[
\delta = \min(1,\gamma) 
\]
then 
 \beq         \label{Dyson:phi-formula-r}
\phi_r(n)  \geq \binom{h-1}{r-1} \delta n
\eeq
for all $r \in \{ 1,\ldots, h\}$. 
\et

This immediately implies Theorem~\ref{Dyson:theorem:Nathanson}. 
 
Like the original proofs of Mann and Artin-Scherk, Dyson's proof is elementary but difficult.  
An important part of Dyson's paper was the introduction of what we call the Dyson transform.

\section{The Dyson transform}

Let $(A_1,\ldots, A_h)$ be an $h$-tuple of sets of positive integers and 
let $n$ be a positive integer with $A_h(n) \geq 1$.  
A \emph{Dyson triple}\index{Dyson triple} is a triple $(a,\ell,c)$ 
such that 
\benum 
\item[(i)]
$c \in A_h$ \text {with } $c \leq n$,
\item[(ii)]
$a \in A_{\ell} \cup \{0\}$ \text{ for some } $\ell \in \{1,\ldots, n-1\}$ or $a > n$,
\item[(iii)]
if $c+a \leq n$, then $c+a \notin A_{\ell}$.
\eenum

Let $(a_0,\ell_0,c_0)$ be a Dyson triple.  We define the set 
\beq                      \label{Dyson:T}
T = T(a_0,\ell_0) = \{c \in A_h: \text{$c \leq n$ and $(a_0,\ell_0,c) $ is a Dyson triple}\}.
\eeq
Note that $c_0 \in T$ and so $T$ is a nonempty subset of $A_h \cap \{1,\ldots, n\}$.   
If $c \in T$ and $c\leq n-a_0$, then 
$c+a_0 \leq n$ and so $c+a_0 \notin A_{\ell_0}$.

The   \emph{Dyson transform}\index{Dyson transform} 
of the $h$-tuple $(A_1,\ldots, A_h)$ with respect to 
the Dyson triple $(a_0,\ell_0,c_0)$ 
is the $h$-tuple $(A'_1,\ldots, A'_h)$ defined as follows: 
\begin{align}
 A'_h &  = A_h \setminus T  \label{Dyson:A'n} \\
A'_{\ell_0} & = A_{\ell_0} \cup (T+a_0)   \label{Dyson:A'ell}  \\
 A'_{\ell}  & = A_{\ell} \qquad   
 \text{for all $\ell \in \{1,\ldots, h-1\} \setminus \{\ell_0\}$.}        \label{Dyson:A'i}  
\end{align}
We have 
\[
 A'_{\ell_0}(m) = A_{\ell_0}(m) + T(m-a_0)  
 \]
for all $m \in \{1,\ldots, n\}$.  If $T(n-a_0) \geq 1$, then $A'_{\ell_0}(n) > A_{\ell_0}(n)$. 
For mathematical induction, the important inequality is 
\beq                    \label{Dyson:induction-An}
A'_h(n) = A_h(n) - T(n)  < A_h(n).  
\eeq 

\bl               \label{Dyson:lemma:0} 
Let $(A'_1,\ldots, A'_h)$ be the Dyson transform of $(A_1,\ldots, A_h)$ with respect to 
the Dyson triple $(a_0,\ell_0,c_0)$.  
Let $m \in \{1,\ldots, n\}$.  
For every set $S$ of positive integers, 
if $p \in \{1,\ldots, m\}$ and $p \in S+A'_{\ell_0}+A'_h$, then $p \in S+A_{\ell_0}+A_h$ and 
\beq              \label{Dyson:0-a} 
(S+A'_{\ell_0}+A'_h)(m) \leq (S +A_{\ell_0}+A_h)(m).
\eeq
\el

\begin{proof}
Let $p \in \{1,\ldots, n\}$.  
We first prove that if $p \in A'_{\ell_0} + A'_h$, 
then $p \in A_{\ell_0} + A_h$.  

If $p \in A'_{\ell_0} + A'_h$, then   
\[
p = x+y
\]
for some  
\[
x \in A'_{\ell_0}   \cup  \{0\} = A_{\ell_0}  \cup  \{0\}  \cup (T+a_0) 
\] 
and 
\[
y \in A'_h   \cup  \{0\}.
\]
If $x \in A_{\ell_0} \cup \{ 0\}$, then 
\[
p = x+y \in A_{\ell_0}+A'_h \subseteq A_{\ell_0}+A_h.
\]
If $x \in T+a_0$, then  $x = c+a_0$ for some  $c \in T \subseteq A_h$.  
The inequality $c > 0$ implies  $a_0 < x$.
If  $y = 0$, then  
\[
p = x = c+a_0 \in A_{\ell_0}+A_h. 
\]
If $y \in A'_h = A_h \setminus T$, then   
\[
1 \leq y \leq   a_0 + y < x+y = p \leq n. 
\]
Because $a_0 \in A_{\ell_0}$ and $y \notin T$, the definition of the set $T$  
implies $a_0 + y \in A_{\ell_0}$ and so 
\[
p = x+y = c+a_0+y \in A_{\ell_0}+A_h.
\]
Thus, if $p \in \{1,\ldots, n\}$ and $p \in A'_{\ell_0} + A'_h$, 
then $p \in A_{\ell_0} + A_h$.

Let  $S$ be a set of positive integers and let $m \in \{1,\ldots, n\}$.  
If $m \in S+A'_{\ell_0}+A'_h$, then $m = s+p$ 
for some $s \in S \cup \{0\}$ and so $p = m-s \in A'_{\ell_0}+A'_h$ 
with $p  \leq  m \leq n$.
It follows that $p \in A_{\ell_0}+A_h$ 
and so $m = s+p \in S+ A_{\ell_0}+A_h$.  Inequality~\eqref{Dyson:0-a} follows immediately. 
This completes the proof.  
\end{proof}

\bc
Let
\[
\phi_h(m) = (A_1+\cdots + A_h)(m) \qqand \phi'_h(m) = (A'_1+\cdots + A'_h)(m). 
\]
Then 
\beq              \label{Dyson:0-b} 
\phi'_h(m)  \leq \phi_h(m) 
\eeq
for all $m \in \{1,\ldots, n\}$. 
\ec

\begin{proof}
Recall that $A_{\ell} = A'_{\ell}$ for $\ell \in \{1,\ldots,n-1\}\setminus \{\ell_0\}$.  
We obtain~\eqref{Dyson:0-b} from~\eqref{Dyson:0-a} 
by choosing 
\[
S = \sum_{\substack{\ell=1 \\ \ell \neq \ell_0}}^{n-1} A_{\ell}  
= \sum_{\substack{\ell=1 \\ \ell \neq \ell_0}}^{n-1} A'_{\ell}.
\]   
This completes the proof. 
\end{proof}

\bl               \label{Dyson:lemma:1} 
Let $(A'_1,\ldots, A'_h)$ be the Dyson transform of $(A_1,\ldots, A_h)$ with respect to 
the Dyson triple $(a_0,\ell_0,c_0)$.  
Let $S$ be a set of positive integers and let $m \in \{1,\ldots, n\}$.  If
\beq          \label{Dyson:up}
m \in \left(  S+A'_{\ell_0} \right) \setminus  \left(  S+A_{\ell_0} \right)
\eeq 
then 
\beq          \label{Dyson:down}
m-a_0 \in  \left(  S+A_h \right) \setminus  \left(  S+A'_h \right). 
\eeq
Moreover,    
\beq          \label{Dyson:updown}
 (S +A'_{\ell_0})(m) + (S +A'_h)(m) \leq  (S  +A_{\ell_0})(m) + (S +A_h)(m).
\eeq
\el

\begin{proof}
We have 
\[
m \in \left( S+A'_{\ell_0} \right) \setminus \left(  S+A_{\ell_0} \right) 
\subseteq  S+  \left( A'_{\ell_0} \setminus  A_{\ell_0} \right) 
\subseteq  S+  T+a_0.  
\]
There exists $s \in S \cup \{0\}$ and $c \in T$ such that 
\[
m = s + c+a_0   
\]
and so 
\[
m-a_0 = s + c \in S + T \subseteq S + A_h.
\] 

If $m - a_0 \in S + A'_h$, then 
\[
m - a_0 = s' + c' 
\]
for some $s' \in S \cup \{0\}$ and $c' \in A'_h \cup \{0\}$.  
If $c' = 0$, then 
\[
m = s' + a_0 \in S + A_{\ell_0} 
\]
which is absurd.  
If $c' \neq 0$, then 
$c' \in A'_h = A_h \setminus T$ and 
\[
1 \leq c' + a_0 = m-s' \leq m \leq n.
\] 
Because $c' \notin T$, we have 
\[
c' + a_0 \in A_{\ell_0} 
\]
and   
\[
m = s' + c' + a_0  \in S + A_{\ell_0}
\]
which is  also absurd. 
Therefore, $ m - a_0 \notin S+A'_h$ and~\eqref{Dyson:up} implies~\eqref{Dyson:down}.

Let $p \in \{1,\ldots, m\}$.  If  
$p \in \left( S+A'_{\ell_0} \right) \setminus \left(  S+A_{\ell_0} \right)$, 
then $p > a_0$ and $p-a_0 \in \{1,\ldots, m - a_0\}$ and 
$p-a_0 \in  \left(  S+A_h \right) \setminus  \left(  S+A'_h \right)$.
It follows that 
\begin{align*}
(S+A'_{\ell_0})(m) - (S+A_{\ell_0})(m) 
& = \left(  \left( S+A'_{\ell_0} \right) \setminus \left(  S+A_{\ell_0} \right) \right)(m) \\
& \leq \left(  \left(  S+A_h \right) \setminus  \left(  S+A'_h \right) \right)(m-a_0) \\
& \leq \left(  \left(  S+A_h \right) \setminus  \left(  S+A'_h \right) \right)(m) \\
& =  (S+A_h)(m) - (S+A'_h)(m)
\end{align*}
which is equivalent to~\eqref{Dyson:updown}. 
This completes the proof. 
\end{proof}

\bt               \label{Dyson:theorem:3zz} 
Let $(A_1,\ldots,A_h)$ be an $h$-tuple of sets of positive integers with $A_h(n) \geq 1$   
and let $(a_0,\ell_0,c_0)$ be a Dyson triple for $(A_1,\ldots,A_h)$.  
Let $(A'_1,\ldots,A'_h)$ be the Dyson transform of $(A_1,\ldots,A_h)$ 
with respect to $(a_0,\ell_0,c_0)$.   For all $m \in \{1,\ldots, n\}$  and $r \in \{ 1,\ldots, h\}$, let 
\[
\phi_r(m) = \sum_{I\in \mci(r,h)} S_I(m) 
\] 
be the rank $r$ counting function for the $h$-tuple of sets $(A_1,\ldots, A_h)$ and let 
\[
\phi'_r(m) = \sum_{I\in \mci(r,h)} S'_I(m) 
\] 
be the rank $r$ counting function for the $h$-tuple of sets $(A'_1,\ldots, A'_h)$.  
Then
\beq                          \label{Dyson:inequality-r}
\phi'_r(m) \leq \phi_r(m).  
\eeq
\et

\begin{proof}
Let $\mci'(r)$ be the set 
of all subsets of size $r$ of the set 
\[
\{1,\ldots, \ell_0 - 1, \ell_0 +1,\ldots, n-1\}. 
\]
For all positive integers $m$, we have 
\begin{align}            \label{Dyson:phi-formula-1}
\phi_r(m) = &  \sum_{I \in \mci'(r)} S_I(m) 
+ \sum_{I \in \mci'(r-1)} \left( (S_I +A_{\ell_0})(m) + (S_I+A_h)(m) \right) \\
& + \sum_{I \in \mci'(r-2)} (S_I +A_{\ell_0} + A_h)(m).        \nonumber
\end{align}
Similarly, 
\begin{align}            \label{Dyson:phi-formula-2}
\phi'_r(m) = &  \sum_{I \in \mci'(r)} S_I(m) 
+ \sum_{I \in \mci'(r-1)} \left( (S_I +A'_{\ell_0})(m) + (S_I+A'_h)(m) \right) \\
& + \sum_{I \in \mci'(r-2)} (S_I +A'_{\ell_0} + A'_h)(m).        \nonumber
\end{align}

Let $m \in \{1,\ldots, n\}$.  
From Lemma~\ref{Dyson:lemma:0} we have 
\beq          \label{Dyson:phi-formula-3}
\left( S_I +A'_{\ell_0} + A'_h \right)(m) \leq \left( S_I +A_{\ell_0} + A_h\right) (m)
\eeq
for all $I \in \mci'(r-2)$. 
From Lemma~\ref{Dyson:lemma:1} we have 
\beq          \label{Dyson:phi-formula-4}
 (S_I +A'_{\ell_0})(m) + (S_I+A'_h)(m) \leq  (S_I +A_{\ell_0})(m) + (S_I+A_h)(m)
\eeq
for all $I \in \mci'(r-1)$.

Equations~\eqref{Dyson:phi-formula-1} 
and~\eqref{Dyson:phi-formula-2} and 
inequalities~\eqref{Dyson:phi-formula-3} and~\eqref{Dyson:phi-formula-4} 
imply
\[
\phi_r'(m) \leq \phi_r(m). 
\]
This completes the proof of Theorem~\ref{Dyson:theorem:3zz}. 
\end{proof}

\section{Minimal Dyson triples}
By inequality~\eqref{Dyson:inequality-r} of Theorem~\ref{Dyson:theorem:3zz}, 
for all $r \in \{1,\ldots, h\}$, 
a lower bound for the rank $r$ counting function 
$\phi'_r(n)$ for the Dyson transform $h$-tuple $(A'_1,\ldots,A'_h)$  
implies a  lower bound for the rank $r$ counting function $\phi_r(n)$ 
for the original $h$-tuple $(A_1,\ldots,A_h)$.  
The next step (Theorem~\ref{Dyson:theorem:4zz}) 
is to prove that a lower bound for the rank $1$ counting 
function $\phi_1(n)$ for the original $h$-tuple $(A_1,\ldots,A_h)$
implies a lower bound for the rank $1$ counting function 
$\phi'_1(n)$ for the Dyson transform $h$-tuple $(A'_1,\ldots,A'_h)$.  
For this, we introduce minimal Dyson triples.    
The core of the construction is the following simple fact.

\bl                                  \label{Dyson:lemma:simple-1}
Let $c > 0$ and let $B$ be a nonempty set of real numbers that is bounded above. 
There exists $b \in B$ such that $b+c \notin B$.
\el

\begin{proof}
Let $b_0 \in B$.    
If $b+c \in B$ for all $b \in B$, then $b_0 + ic \in B$ for all integers $i \geq 0$ 
and so $b_0+ ic \leq \sup B$ for all $i \geq 0$, which is absurd.  
This completes the proof. 
\end{proof}

\bl                                     \label{Dyson:lemma:simple-2}
Let $(A_1,\ldots, A_h)$ be an $h$-tuple of sets of positive integers. 
Let $n$ be a positive integer and let $c \in A_h$ with $c \leq n$.  
For all $\ell \in \{1, 2, \ldots, h-1\}$, there exists $a_{\ell}^{\sharp} \in A_{\ell} \cup \{ 0\}$ 
such that $a_{\ell}^{\sharp} \leq n$ and $(a_{\ell}^{\sharp},\ell,c)$ is a Dyson triple.  
\el

\begin{proof}
Apply Lemma~\ref{Dyson:lemma:simple-1} to the set $B = \{ a \in A_{\ell} \cup \{ 0\}: a \leq n\}$. 
\end{proof}

Let $a_0$ be the smallest nonnegative integer such that there exists 
a Dyson triple $(a,\ell,c)$ with $a = a_0$. 
A \emph{minimal Dyson triple} 
is a Dyson triple of the form $(a_0,\ell_0,c_0)$ with $a_0 \in A_{\ell_0} \cup \{0\}$. 
The integer $a_0$ is unique 
but minimal Dyson triples are not unique.

\bl               \label{Dyson:lemma:simple-4}
Let $n$ be a positive integer with $A_h(n) \geq 1$ and let  
$(a_0,\ell_0,c_0)$ be a minimal Dyson triple.  Then 
\beq                 \label{Dyson:a0-bound} 
0 \leq a_0 \leq n.
\eeq
\el

\begin{proof}
This follows immediately from Lemma~\ref{Dyson:lemma:simple-2}. 
\end{proof}

\bl                       \label{Dyson:lemma:mc}
Let $(A_1,\ldots, A_h)$ be an $h$-tuple of sets of positive integers with $A_h(n) \geq 1$  
and let $(a_0,\ell_0,c_0)$ be a minimal Dyson triple. 
If $m$ and $c$ are  integers such that 
\beq                 \label{Dyson:p} 
c \in A_h \qqand m - a_0 < c \leq m \leq n
\eeq
then 
\beq                 \label{Dyson:q} 
A_{\ell}(m) \geq A_{\ell}(m-c) +1 +A_{\ell}(c-1) 
\eeq 
for all $\ell \in \{1,\ldots, n-1\}$. 
\el

\begin{proof}
Let $m$ and $c$ be integers that satisfy~\eqref{Dyson:p}.  
Then 
\[
m-c <  a_0.
\]  
For all $\ell \in \{1,\ldots, n-1\}$, there are $A_{\ell}(m-c)+1$ integers 
$a \in A_{\ell} \cup \{0\}$ such that 
\[
0 \leq a \leq m-c. 
\]
For each such $a$ we have    
\[
a < a_0
\]
and 
\[
c \leq a+c \leq m \leq n. 
\]  
The minimality of $a_0$ in the Dyson triple $(a_0,\ell_0,c_0)$ implies that 
$(a,\ell,c)$ is not a Dyson triple and so $a+c \in A_{\ell}$.  
Therefore,
\[
A_{\ell}(m)-A_{\ell}(c-1) \geq A_{\ell}(m-c)+1.
\]
This completes the proof. 
\end{proof}


\bt                       \label{Dyson:theorem:4zz}   
Let $(A_1,\ldots,A_h)$ be an $h$-tuple of sets of positive integers with $A_h(n) \geq 1$   
and let $(a_0,\ell_0,c_0)$ be a minimal Dyson triple for $(A_1,\ldots,A_h)$.
Let $(A'_1,\ldots,A'_h)$ be the Dyson transform of $(A_1,\ldots,A_h)$ 
with respect to $(a_0,\ell_0,c_0)$.  
Consider the rank 1 counting functions  
\[
\phi_1(m) = \sum_{i=1}^n A_i(m) 
\] 
and  
\[
\phi'_1(m) =  \sum_{i=1}^n A'_i(m). 
\] 
If 
\[
\phi_1(m) \geq \gamma m  
\]
for all $m \in \{1,\ldots, n\}$ and 
\[
\delta = \min(1,\gamma) 
\] 
then 
\[
\phi'_1(m) \geq \delta m 
\] 
for all $m \in \{1,\ldots, n\}$.
\et

\begin{proof} 
Recall that 
\[
T = \{c\in A_h: \text{$c\leq n$ and $(a_0,\ell_0,c)$ is a Dyson triple} \}
\]
and 
\[
A'_h = A_h \setminus T \qqand A'_{\ell_0} = A_{\ell_0} \cup (T+a_0). 
\] 
Let $m \in \{1,\ldots, n\}$.  Because $T \subseteq A_h$, we have 
\[
A'_h(m) = A_h(m) - T(m).  
\]
Moreover, $0 \leq a_0 \leq n$ by Lemma~\ref{Dyson:lemma:simple-4}. 
 
If $T(m-a_0) = 0$, then $(T+a_0)(m) = 0$ and $A'_{\ell_0}(m) = A_{\ell_0}(m)$.  
It follows that 
\begin{align*}
\phi'_1(m) 
& = \sum_{\ell=1}^{n-1} A'_{\ell}(m) + A'_h(m)  \\
& =  \sum_{\ell=1}^{n-1} A_{\ell}(m) + A_h(m) - T(m)  \\
& = \phi_1(m) - T(m) \\
& \leq \phi_1(m).
\end{align*}

If $T(m-a_0) > 0$, then there exist integers $c \in T$ such that 
$c+a_0 \leq m \leq n$ and so  $c + a_0 \notin A_{\ell_0}$ and $c+a_0 \in A'_{\ell_0}$. 
Thus, 
\[
A'_{\ell_0}(m) - A_{\ell_0}(m) \geq T( m - a_0).
\]  
Conversely, if $x \in A'_{\ell_0} \setminus A_{\ell_0}$ and $x \leq m$, 
then $x \in T+a_0$ and so  
$x = c + a_0$ for some $c \in T$ with $c \leq m - a_0$.  Thus, 
\[
A'_{\ell_0}(m) - A_{\ell_0}(m) \leq T( m - a_0).
\]  
It follows that 
\[
A'_{\ell_0}(m) = A_{\ell_0}(m) + T( m - a_0).
\] 
The inequality 
\[
T(m) \geq T(m-a_0) 
\]
implies 
\begin{align*}
\phi'_1(m) 
& = \sum_{\substack{\ell=1\\\ell \neq \ell_0}}^{n-1} A'_{\ell}(m) + A'_h(m) + A'_{\ell_0}(m)  \\
& = \sum_{\substack{\ell=1\\\ell \neq \ell_0}}^{n-1} A_{\ell}(m) + \left( A_h(m) - T(m)  \right) 
+ \left( A_{\ell_0}(m) + T( m - a_0)  \right) \\  
& =  \sum_{\ell =1}^n A_{\ell}(m)  -  T(m) + T(m-a_0)  \\
& = \phi_1(m)  - T(m) + T(m-a_0)  \\
&  \leq \phi_1(m). 
\end{align*}
If $T(m) = T(m - a_0)$, then 
\beq           \label{Dyson:lowerLower}
\phi'_1(m) = \phi_1(m) \geq \gamma m \geq \delta m.
\eeq 
If 
\[
T(m) > T(m - a_0)   
\]
then $\phi'_1(m) < \phi_1(m)$ and we cannot claim that  $ \phi'_1(m) \geq \gamma m$.  
However, we can prove that $\phi'_1(m) \geq \delta m$.

Because $T$ is a subset of $A_h$, we have 
\[
0 < T(m) -T(m-a_0) \leq  A_h(m) - A_h(m-a_0).  
\]
There is a smallest integer $b$ in $A_h$ such that 
\beq       \label{Dyson:b-minimal}
m-a_0 < b \leq m.
\eeq
Then $A_h(b-1) = A_h(m-a_0)$ and  
\begin{align*}
\phi'_1(m) & = \phi_1(m)  - T(m) +T(m-a_0) \\
& \geq  \phi_1(m)  - A_h(m) +A_h(m-a_0) \\
& =  \phi_1(m)  - A_h(m) +A_h(b-1).
\end{align*}
Applying Lemma~\ref{Dyson:lemma:mc} with $c = b$, we obtain 
\[
A_{\ell}(m) \geq A_{\ell}(m-b) +1 +A_{\ell}(b-1) 
\]  
for all $\ell \in \{1,\ldots, n-1\}$ and so 
\begin{align*}
\phi'_1(m) & \geq \phi_1(m)  - A_h(m) +A_h(b-1) \\
& = \sum_{\ell =1}^{h-1} A_{\ell}(m) +A_h(b-1) \\ 
& \geq \sum_{\ell=1}^{h-1} \left( A_{\ell}(m-b) +1 +A_{\ell}(b-1) \right)  +A_h(b-1) \\
& = \phi_1(b-1) + \phi_1(m-b)  - A_h(m-b) + h-1.
\end{align*} 
From the inequality $\phi_1(m) \geq \gamma m$ for all $m \in  \{ 1,\ldots, n \}$,  
we obtain  
\beq         \label{Dyson:r}
\phi'_1(m) \geq \gamma (b-1)  + \phi_1(m-b) - A_h(m-b) + h-1. 
\eeq

We shall prove that there exists no positive integer $z$ such that 
\beq                                            \label{Dyson:s}
z < a_0  \qqand \phi_1(z) - A_h(z) < \delta z.
\eeq
If such integers exist, let $z$ be the smallest one.  The inequality  $\delta \leq \gamma$ implies 
\beq                  \label{Dyson:v}
A_h(z) > \phi_1(z) - \delta z \geq \gamma z - \delta z = (\gamma - \delta)z \geq 0.
\eeq
Because $A_h(z) > 0$, there exists $c \in A_h$ with $c \leq z$.  
We have
\[
z-a_0 < 0 < c \leq  z \leq n 
\]
and so $c$ and $z$ satisfy inequality~\eqref{Dyson:p}. 
Applying Lemma~\ref{Dyson:lemma:mc}, we obtain 
\begin{align}                  \label{Dyson:t}
\phi_1(z) - A_h(z) 
& = \sum_{\ell =1}^{h-1} A_{\ell}(z)     \nonumber \\ 
&  \geq  \sum_{\ell =1}^{h-1} \left(  A_{\ell}(z-c) +1 +A_{\ell}(c-1)  \right)       \nonumber \\ 
& = \phi_1(z-c) - A_h(z-c) + \phi_1(c-1) - A_h(c-1) + h-1.
\end{align}
The nonnegative integers $z-c$ and $c-1$ are less than $z$, which is the smallest  
positive integer satisfying~\eqref{Dyson:s}.       
It follows that, if $z-c$ and $c-1$ are positive, then 
\[
\phi_1(z-c) - A_h(z-c) \geq \delta (z-c)
\]
and 
\[
\phi_1(c-1) - A_h(c-1) \geq \delta (c-1). 
\]
These inequalities also hold if $z-c=0$ or if $c-1=0$.
Inserting the inequalities into~\eqref{Dyson:t} gives 
\begin{align}            
\phi_1(z) - A_h(z) & \geq \delta (z-c)+ \delta (c-1)+ h-1 \nonumber \\
& = \delta (z-1)+ h-1.                  \label{Dyson:u}
\end{align}                 
Because $h \geq 2$ and $\delta \leq 1$, 
inequalities~\eqref{Dyson:v} and~\eqref{Dyson:u}
imply 
\[
\delta z > \phi_1(z) - A_h(z) \geq \delta z+ h-1- \delta \geq  \delta z
\]
which is absurd.  
Therefore, no integer $z$ satisfies inequalities~\eqref{Dyson:s}.  

Recall that $b$ is the smallest integer in the set $A_h$ such that 
\[
m - a_0 < b \leq m 
\]
and so $m-b < a_0$.  The impossibility of~\eqref{Dyson:s} implies 
\[
\phi_1(m-b) - A_h(m-b) \geq \delta (m-b). 
\]
Because $h \geq 2$ and $\delta = \min(1,\gamma)$, inequality~\eqref{Dyson:r} gives 
\begin{align*}
\phi'_1(m) 
& \geq \gamma (b-1)  + \phi_1(m-b) - A_h(m-b) + h-1 \\
& \geq \gamma (b-1)  +\delta (m-b) + h-1 \\
& = \delta m + (\gamma - \delta)(b-1) + h-1-\delta \\
& \geq \delta m  + h-1-\delta \\
& \geq \delta m.
\end{align*} 
This completes the proof of Theorem~\ref{Dyson:theorem:4zz}. 
\end{proof}

\section{Completion of the proof of Dyson's theorem}

Dyson's proof of Theorem~\ref{Dyson:theorem:Dyson} 
is by double induction on $h$ and $A_h(n)$.    
The theorem is trivially true for $h=1$.  
Let $h \geq 2$ and assume Theorem~\ref{Dyson:theorem:Dyson} 
is true for every $(h-1)$-tuple of sets of positive integers.
For all $r \in \{1,\ldots, h-1\}$, 
let $\widetilde{\phi}_r(m)$ be the 
rank $r$ counting function for the $(h-1)$-tuple of sets $(A_1,\ldots, A_{h-1})$. 
For all $I =  \{i_1,\ldots, i_{r-1}\} \in \mci(r-1,h-1)$, let 
\[
\widetilde{I} =  \{i_1,\ldots, i_{r-1},h\}  \in \mci(r-1,h) 
\]
and 
\[
\widetilde{\mci}(r,h) = \left\{ \widetilde{I} : I \in \mci(r-1,h-1) \right\}.
\]
Then
\[
\mci(r,h) = \mci(r,h-1) \cup \widetilde{\mci}(r,h)
\]
and 
\[
\mci(r,h-1) \cap \widetilde{\mci}(r,h) = \emptyset.
\]

We proceed by induction on $A_h(n)$.  
If  $A_h(n) = 0$, then, for all $I = (i_1,\ldots, i_{r-1})\in \mci(r-1,h-1)$ 
and $m \in \{1,\ldots, n\}$, we have 
\begin{align*}
S_I(m)  & = \left( A_{i_1}+\cdots + A_{i_{r-1}} \right)(m) \\
& = \left( A_{i_1}+\cdots + A_{i_{r-1}} + A_h \right)(m) \\
& = S_{\widetilde{I}}(m)  
\end{align*} 
and so 
\[
\sum_{\widetilde{I} \in \widetilde{\mci}(r,h)}S_{\widetilde{I}}(m) = \sum_{I\in \mci(r-1,h-1)}S_I(m) 
=   \widetilde{\phi}_{r-1}(m). 
\]
Applying the induction hypothesis for $h-1$ sets and $r \in \{1,\ldots, h-1\}$, we obtain 
\begin{align*}
\phi_r(m) & = \sum_{I\in \mci(r,h)}S_I(m) \\ 
& = \sum_{I\in \mci(r,h-1)}S_I(m) + \sum_{\widetilde{I} \in \widetilde{\mci}(r,h)}S_{\widetilde{I}}(m) \\
& =  \widetilde{\phi}_{r}(m) +  \widetilde{\phi}_{r-1}(m)\\
& \geq \left( \binom{h-2}{r-1} + \binom{h-2}{r-2} \right) \delta m \\
& = \binom{h-1}{r-1} \delta m. 
\end{align*}
For $r=h$ we have 
\begin{align*}
\phi_h(m) & = (A_1+\cdots + A_{h-1} + A_h)(m) \\ 
& = (A_1+\cdots + A_{h-1})(m) 
= \widetilde{\phi}_{h-1}(m) \\
& \geq \delta m.
\end{align*} 
This proves Theorem~\ref{Dyson:theorem:Dyson} when $A_h(n)=0$.

Suppose that $A_h(n) \geq 1$.  
Let $(a_0,\ell_0,c_0)$ be a minimal Dyson triple for the $h$-tuple $(A_1,\ldots, A_h)$ 
and let  $(A'_1,\ldots, A'_h)$ be the Dyson transform of 
$(A_1,\ldots, A_h)$ with respect to $(a_0,\ell_0,c_0)$.  
For all $r \in \{1,\ldots, h\}$, let $\phi'_r(m)$ be the rank $r$ counting function 
of the $h$-tuple $(A'_1,\ldots, A'_h)$.  

We have $\phi_1(m) \geq \gamma m$ for all $m \in \{1,\ldots, n\}$. 
Let $\delta = \min(1,\gamma)$. 
By Theorem~\ref{Dyson:theorem:4zz}, 
\[
\phi'_1(m) \geq \delta m
\]
for all $m \in \{1,\ldots, n\}$. 
By inequality~\eqref{Dyson:induction-An},
\[
A'_h(n) < A_h(n).  
\]
Let $m \in \{1,\ldots, n\}$.  Applying the induction hypothesis to the $h$-tuple $(A'_1,\ldots, A'_h)$,   
we obtain 
\[
\phi'_r(m) \geq \delta m. 
\] 
By Theorem~\ref{Dyson:theorem:3zz}, 
\[
\phi_r(m) \geq \phi'_r(m).
\]
This completes the proof.

\section{Asymptotic density and Kneser's theorem} 
Let $A$ and $B$ be sets of positive integers.  
By Mann's $\alpha +\beta$ theorem, Shnirel'man density is superadditive 
in the following sense:  
\[
\sigma(A+B) \geq \min(1,\sigma(A) + \sigma(B)).
\]
The \emph{lower asymptotic density}\index{lower asymptotic density}\index{asymptotic density!lower} 
of the set $A$ is 
\[
d_L(A) = \liminf_{n = 1,2,\ldots} \frac{A(n)}{n}. 
\]
Asymptotic density is not superadditive.  For example, let $k$, $\ell$, and $m$ be positive integers
 such that 
 \[
 k+\ell \leq m.  
 \]
Consider the sets 
\[
A = \left\{ n \in \N: n \equiv r \pmod{m} \text{ for some } r \in \{0,1,\ldots, k-1\} \right\}
\]
and 
\[
B = \left\{ n \in \N: n \equiv s \pmod{m} \text{ for some } s \in \{0,1,\ldots, \ell - 1\} \right\} 
\]
with lower asymptotic densities 
\[
d_L(A) = \frac{k}{m} \qqand d_L(B) = \frac{\ell}{m}.  
\]
The sumset 
\[
A+B = \left\{ n \in \N: n \equiv t \pmod{m} \text{ for some } t \in \{0,1,\ldots, k+\ell - 2\} \right\}.
\]
has lower asymptotic density 
\begin{align*}
d_L(A+B) & = \frac{k+\ell-1}{m} \\
& = d_L(A) + d_L(B) - \frac{1}{m} \\
& <  d_L(A) + d_L(B).
\end{align*}

In 1953, Martin Kneser~\cite{knes53} proved the deep and beautiful theorem that 
the only counterexamples to superadditivity for lower asymptotic density of sumsets 
are pairs of sets $A$ and $B$ that are essentially unions of congruence classes modulo $m$
for some positive integer $m$.  The long and difficult proof is a \emph{tour de force} based 
on the masterful application of iterated Dyson  transforms. 
For proofs of Kneser's theorem, see Kneser~\cite{knes53} and also 
Halberstam-Roth~\cite{halb-roth66}, Mann~\cite{mann65}, and Nathanson~\cite{nath25aa}.

\section{The  $e$-transform in an  abelian group}
Analogous to his theorem on the lower asymptotic density 
of sums of sets of integers, Kneser also proved a theorem 
on sums of finite subsets of an additive abelian group. 
The proof uses a simple form of the Dyson transform.  
Let $G$ be an abelian group, written additively, with identity element $0$.
For $X \subseteq G$ and $e \in G$, we define the \emph{translates}\index{translate} 
\[
X+e = \{x+e:x \in X\} 
\] 
and
\[
X - e = \{x - e:x \in X\}.
\] 
 Let $(A,B)$ be a pair of nonempty subsets of $G$.  
For all $e \in A$, let
\[
A(e) = A \cup (B+e)
\]
and 
\[
B(e) = B \cap (A-e)
\]
The pair of sets $(A(e),B(e))$ is the \emph{$e$-transform}\index{Dyson $e$-transform} 
of the pair $(A,B)$.  
Note that  
\[
A \subseteq A(e) 
\]
and 
\[
 B(e) \subseteq B.
\]


\bt                            \label{Dyson:theorem:e-transform}   
Let $A$ and $B$ be nonempty subsets of the abelian group $G$ and let $e \in G$.  Then 
\beq                   \label{Dyson:e-transform-1}   
A(e) + B(e) \subseteq A+B 
\eeq
and
\beq                   \label{Dyson:e-transform-2}   
A(e) \setminus A = e + \left(  B\setminus B(e) \right).
\eeq
If $A$ and $B$ are finite sets, then 
\beq                   \label{Dyson:e-transform-3}   
|A| + |B| = |A(e)| + |B(e)|.
\eeq If $e \in A$ and $0 \in B$, then $e \in A(e)$ and $0 \in B(e)$. 
\et

\begin{proof}
Let $a'+b' \in A(e)+B(e)$, where $a' \in A(e)$ and $b' \in B(e)$. 
Then $a' \in A$ or $a' \in B+e$.
If $a' \in A$, then $b' \in B(e) \subseteq B$ implies $a'+b' \in A+B$.  
If $a' \in B+e$, then $a' = b'' + e$ for some $b'' \in B$.
We have $b' \in B(e)$ and so $b' = a'' - e$ for some $a'' \in A$.
It follows that 
\[
a'+b' = (b''+e)+(a''-e) = a'' + b'' \in A+B.
\]
This proves~\eqref{Dyson:e-transform-1}. 

We have 
\[
b \in  B\setminus B(e) 
\]
if and only if 
\[
b \in B \qqand b \notin A-e 
\]
if and only if 
\[
b+e \in B +e\qqand b +e \notin A 
\]
if and only if 
\[
b+e \in A(e) \setminus A. 
\]
This proves~\eqref{Dyson:e-transform-2}.

If $A$ and $B$ are finite sets, then~\eqref{Dyson:e-transform-2} implies  
\[
\left|  A(e)   \right|  -  \left| A  \right| 
= \left|  A(e) \setminus A  \right| = \left|  B \setminus  B(e) \right|
= \left|  B  \right|  -  \left|  B(e) \right|.
\]
This proves~\eqref{Dyson:e-transform-3}. 

 If $e \in A$ and $0 \in B$, then $e \in A(e)$ and $0 \in B \cap (A-e) = B(e)$.  
 This completes the proof. 
\end{proof}

As an example of the use the $e$-transform in additive number theory, 
we prove I. Chowla's generalization of the Cauchy-Davenport theorem.  
We begin with a simple counting argument.

\bl                \label{Dyson:lemma:>|G|}
Let $G$ be a finite group, written multiplicatively and not necessarily abelian.  
Let $A$ and $B$ be subsets of $G$ and let 
\[
AB = \{ab:a\in A \text{ and } b \in B\}.
\]
If 
\[
|A| + |B| > |G|
\]
then 
\[
AB = G.
\]
\el

\begin{proof}
Let $x \in G$ and let 
\[
xB^{-1} = \{xb^{-1}:b \in B\}.
\]
Then 
\[
\left| xB^{-1} \right| = |B|
\]
and  
\[
|A| +\left| xB^{-1} \right| = |A| + |B|  > |G|. 
\]
By the pigeonhole principle, the sets $A$ and $ xB^{-1}$ are not disjoint, 
and so there exist $a \in A$ and $b\in B$ such that $a = xb^{-1}$.
Thus, $x = ab \in AB$ and $AB=G$.
This completes the proof.  
\end{proof}

\bt [I. Chowla]              \label{Dyson:theorem:IChowla}
Let $m \geq 2$ and let $A$ and $B$ be nonempty subsets of 
the additive abelian group $\Z/m\Z$ of congruence classes modulo $m$.  
If $0 \in B$ and $\gcd(b,m) = 1$ for all $b \in B \setminus \{0\}$, then 
\beq                   \label{Dyson:IChowla}  
|A+B| \geq \min(m,|A|+|B|-1).
\eeq
\et

\begin{proof}
 By Lemma~\ref{Dyson:lemma:>|G|}, if $|A|+|B| > m$, then $A+B = \Z/m\Z$ 
 and $|A+B| = m$.  Thus, we can assume that 
 \[
 |A|+|B| \leq m.
 \]
Inequality~\eqref{Dyson:IChowla}  also holds if $|A|=1$ or $|B|=1$, and so we can assume 
that 
\[
|A| \geq 2 \qqand |B| \geq 2.
\]
 
If the theorem is false, then there is a counterexample, that is, a pair of sets $(A,B)$ such that 
$|A| \geq 2$ and $|B| \geq 2$ and 
 \[
|A+B| <  |A|+|B|-1.
 \]
 Choose the counterexample $(A,B)$ so that the cardinality of $B$ is minimal.  
 Because $|B| \geq 2$, there exists $b^* \in B$ with $b^* \neq 0$ and $\gcd(b^*,m)=1$.  
 
If $a+b^* \in A$ for all $a \in A$, then $a + jb^* \in A$ for all $j = 0,1,2,3,\ldots$. 
Because the set $A$ is finite, there exist integers $j_1 < j_2$ with $j_0 = j_2-j_1$ 
minimal such that 
$a + j_1b^* = a + j_2b^*$.  It follows that  
\[
j_0b^* = \left( a + j_1b^*\right)  - \left( a + j_2b^*\right) = 0. 
\]
The divisibility condition $\gcd(b^*,m)=1$ implies $j_0 = m$ and so $A = G$, 
which is absurd.  Therefore, there exists $e \in A$ with $e+b^* \notin A$.  

We apply this element $e$  to the pair $(A,B)$ 
to obtain the $e$-transform   $(A(e),B(e))$. 
By Theorem~\ref{Dyson:theorem:e-transform}, 
\begin{align*}
|A(e)+B(e)| & \leq |A+B| \\ 
& < |A| + |B| -1 \\
& = |A(e)| + |B(e)| - 1
\end{align*}
and also  $0 \in B(e)$.
Moreover, $B(e) \subseteq B$ implies $\gcd(b,m) = 1$ for all $b \in B(e) \setminus \{0\}$. 
Thus, the pair $(A(e),B(e))$ is also a counterexample.  

Because $e + b^* \notin A$, we have $b^* \notin A-e$ and so $b^* \notin B(e)$. 
Thus, $B(e)$ is a proper subset of $B$, which contradicts the minimality of $|B|$.  
This completes the proof. 
\end{proof}

From I. Chowla's theorem we immediately obtain the Cauchy-Davenport theorem. 

\bt [Cauchy-Davenport]              \label{Dyson:theorem:Cauchy-Davenport}
Let $p$ be a prime number and let $A$ and $B$ be nonempty subsets of $\Z/p\Z$.  
Then 
\beq                   \label{Dyson:Cauchy-Davenport}  
|A+B| \geq \min(p,|A|+|B|-1).
\eeq
\et

\begin{proof}
Let  $A$ and $B$ be nonempty subsets of $\Z/p\Z$.  Let $b \in B$ and let $B' = B-b$.  
Then $|B'| = |B|$ and $|A+B'| = |A+B - b| = |A+B|$.  Also, $0 \in B'$.  

Because $p$ is prime, we have $\gcd(b',p) = 1$ for all $ b' \in B'$ with $b' \neq 0$. 
Applying Theorem~\ref{Dyson:theorem:IChowla} with $m = p$ to the sets $A$ and $B'$,   
we obtain 
\[
|A + B| = |A + B'| \geq \min(p,|A|+|B'|-1) = \min(p,|A|+|B|-1).
\]
This completes the proof.  
\end{proof}

Martin Kneser used the $e$-transform to prove the following important group theoretic inequality.

\bt[Kneser]                            \label{Dyson:theorem:Kneser}
Let $G$ be an abelian group, written additively,  and let $A$ and $B$ 
be nonempty finite subsets of $G$.  
If $|A| + |B| \leq |G|$, then there exists a proper subgroup $H$ of $G$ such that 
\[
|A+B| \geq |A| + |B| - |H|.
\]
\et

\begin{proof}
The proof is by induction on $|B|$.  If $|B|=1$, then 
\[
|A+B| = |A| \geq |A| + 1- |H|
\]
for every subgroup $H$ of $G$. 

Let $k \geq 2$ and assume that Kneser's inequality is valid for all pairs $(A,B)$ of 
nonempty finite subsets of $G$ with $|B| < k$. 

Let $(A,B)$ be a pair of nonempty finite subsets of $G$ with $|B| = k$.  
There are two cases.

In the first case, we have 
\[
a+b-b' \in A 
\]
for all $a \in A$ and $b,b'\in B$ and so  
\[
A+b-b' \subseteq A
\]
for all   $b,b'\in B$.  Because $A$ is finite and $|A+b-b'| = |A|$, it follows that 
\[
A+b-b'  =  A
\]
for all   $b,b'\in B$.  
Let $H$ be the subgroup of $G$ generated by the difference set $B-B = \{b-b':b,b' \in B\}$.  
Note that $|B| \leq |B-B| \leq |H|$ and that $A \neq G$ because $|A| + |B| \leq |G|$.  
We have 
\[
A+H = A \neq G
\]
and so $H \neq G$.  Thus, $H$ is a proper subgroup of $G$ such that 
\[
|A+B| \geq |A| \geq |A| + |B| - |H|.
\]

In the second case, there exist $a \in A$ and $b, b' \in B$ such that 
\[
a+b-b' \not \in A .
\]
Let 
\[
e = a - b'.
\] 
We  apply the $e$-transform to the pair $(A,B)$ to obtain the pair $(A(e), B(e))$, where 
\[
A(e) = A \cup (B+e)
\]
and
\[
B(e) = B \cap (A-e).
\]
We have $b+e \in B+e$ and $b+e = a+b-b' \notin  A$ and so $b \notin A - e$.   
Thus, $|A(e)| > |A|$ and $  |B(e)| < |B|$. 
Also,  $b' = a-(a-b') = a - e \in A-e$ and so $b' \in B\cap (A-e) = B(e) \neq \emptyset$. 
Therefore, $B(e)$ is a nonempty finite subset of $G$ with $|B(e)| < |B|$.  
Applying  the induction hypothesis to the pair $(A(e),B(e))$ and applying also 
Theorem~\ref{Dyson:theorem:e-transform}, 
we obtain a  proper subgroup $H$ of $G$ such that 
\[
|A+B| \geq|A(e)+B(e)| 
\geq |A(e)| + |B(e)| - |H| = |A| + |B| - |H|.
\]
This completes the proof. 
\end{proof}

If $p$ is a prime number and $G = \Z/p\Z$, then the only  proper subgroup of $\Z/p\Z$ 
is $H = \{0\}$, and so either $A+B = G$ or $|A+B| \geq |A|+|B| -1$.
Thus, Kneser's inequality immediately implies the Cauchy-Davenport theorem 
(Theorem~\ref{Dyson:theorem:Cauchy-Davenport}).

Kneser also proved a more precise version of Theorem~\ref{Dyson:theorem:Kneser}. 
The \emph{stabilizer}\index{stabilizer} of a nonempty subset $X$ of an abelian group $G$ 
is the set 
\[
H(X) = \{g \in G: g+X = X\}.
\]
The stabilizer of $X$ is a subgroup of $G$.

\bt[Kneser]                           \label{Dyson:theorem:Kneser-H}
Let $G$ be an abelian group, written additively,   let $A$ and $B$ 
be finite nonempty subsets of $G$, and let $H = H(A+B)$ be the stabilizer 
of the sumset $A+B$.   
If 
\[
|A+B| < |A| + |B| 
\]
then 
\[
|A+B| =  |A+H| + |B+H| - |H|.
\]
\et

\begin{proof}
For Kneser's proofs, see~\cite{knes54,knes55} and also Nathanson~\cite[pp. 110--116]{nath96bb}.
\end{proof}

\section{Sums of measurable sets}
A. M. Macbeath proved the first density addition theorem for sums of Lebesgue measurable sets.  
Sierpinski~\cite{sier20} had previously shown that the sum of Lebesgue measurable sets 
is not necessarily Lebesgue measurable.  Introducing inner Lebesgue measure, 
Macbeath obtained a continuous analogue of Mann's theorem 
for the sum of two measurable sets of positive real numbers.  Let $\mu$ denote Lebesgue 
measure and let $\mu_*$ denote inner Lebesgue measure.  

\bt[Macbeath~\cite{macb60a}] 
Let $A$ and $B$ be Lebesgue measurable sets of positive real numbers such that $\inf A = \inf B = 0$ 
and let $A+B=C$. Let  $x$ and $\gamma$ be real numbers such that 
\[
\mu(A \cap [0,t]) + \mu(B \cap [0,t]) \geq \gamma t
\]
for all $t \leq x$. 
 Let $\delta = \min(1,\gamma)$.  Then 
\[
\mu_*(C \cap [0,x]) \geq \delta x. 
\]
\et

Extending this result to $r$-fold sums of $h$ sets of integers, Nathanson obtained a continuous 
analogue of Dyson's theorem.

\bt[Nathanson~\cite{nath80}]
Let $(A_1,\ldots, A_h)$ be an $h$-tuple of Lebesgue measurable sets of positive real numbers 
such that $\inf A_i = 0$ for all $i \in \{1,\ldots, h\}$.  
For all $r \in \{1,\ldots, h\}$, let 
\[
\Phi_r(t) = \sum_{I \in \mci(r,h)} \mu_*(S_I \cap [0,t])
\]
Let  $x$ and $\gamma$ be real numbers such that 
\[
\Phi_1(t) = \mu\left(A_1 \cap [0,t]\right) + \cdots + \mu\left(A_h \cap [0,t]\right) \geq \gamma t
\]
for all $t \leq x$.  Let $\delta = \min(1,\gamma)$.  Then 
\[
\Phi_r(x) \geq \binom{h-1}{r-1} \delta x. 
\]
\et

There is a large literature investigating the Haar measure of sumsets in locally compact 
abelian groups.  Less understood is the density and size of sums of sets of 
nonnegative lattice points in $\Z^n$.


\end{document}